\documentclass[11pt]{article}
\usepackage{graphicx}
\def\be{\begin{equation}}
\def\ee{\end{equation}}
\def\ff#1{\mbox{\boldmath $#1$} }
\def\a{\alpha}
\def\b{\beta}
\def\lam{\lambda}
\def\e{\epsilon}
\def\x{\mbox{\boldmath $x$} }
\newcommand{\mat}[1]{{\left( \begin{array}{cccc}#1\end{array}\!\right)}}

\def\={\approx}

\begin{document}

\title{Metaheuristic Optimization: Algorithm Analysis and Open Problems}

\author{Xin-She Yang \\
Mathematics and Scientific Computing,
National Physical Laboratory,  \\   Teddington, Middlesex TW11 0LW, UK.  }

\date{}
\maketitle

\begin{abstract}
Metaheuristic algorithms are becoming an important part of modern optimization.
A wide range of metaheuristic algorithms
have emerged over the last two decades, and many metaheuristics such as particle swarm optimization
are becoming increasingly popular. Despite their popularity, mathematical analysis of these
algorithms lacks behind. Convergence analysis still remains unsolved for the majority
of metaheuristic algorithms, while efficiency analysis is equally challenging. In this paper, we intend
to provide an overview of convergence and efficiency studies of metaheuristics, and try to provide a
framework for analyzing metaheuristics in terms of convergence and efficiency. This can form
a basis for analyzing other algorithms. We also outline some open questions
as further research topics. \\[10pt]

{\bf Citation Details}:
Yang, X. S., (2011). Metaheuristic optimization: algorithm analysis and open problems, in: Proceedings of 10th International Symposium on Experimental Algorithms (SEA 2011) (Eds. P. M. Pardalos and S. Rebennack), Kolimpari, Chania, Greece, May 5-7 (2011), Lecture Notes in Computer Sciences, Vol. 6630, pp. 21-32.

\end{abstract}


\section{Introduction}

Optimization is an important subject with many important application, and
algorithms for optimization are diverse with a wide range of successful
applications \cite{Floudas,Floudas2}. Among these optimization algorithms, modern metaheuristics
are becoming increasingly popular,
leading to a new branch of optimization, called metaheuristic optimization.
Most metaheuristic algorithms are nature-inspired \cite{Dorigo,Talbi,Yang},
from simulated annealing \cite{Kirk} to ant colony optimization \cite{Dorigo},
and from  particle swarm optimization \cite{Kennedy} to cuckoo search \cite{Yang4}.
Since the appearance of swarm intelligence algorithms such as PSO in the 1990s,
more than a dozen new metaheuristic algorithms have been developed
and these algorithms have been applied to almost all areas of optimization, design,
scheduling and planning, data mining, machine intelligence, and many others.
Thousands of research papers and dozens of books have
been published \cite{Dorigo,Enge,Floudas2,Holland,Talbi,Yang,Yang2}.

Despite the rapid development of metaheuristics, their mathematical
analysis remains partly unsolved, and many open problems need urgent attention.
This difficulty is largely due to the fact the interaction of various components
in metaheuristic algorithms are highly nonlinear, complex, and stochastic.
Studies have attempted to carry out convergence analysis \cite{Auger,Neumann},
and some important results concerning PSO were obtained \cite{Clerc}.
However, for other metaheuristics such as firefly algorithms and ant colony optimization,
it remains an active, challenging topic.
On the other hand, even we have not proved or cannot prove their convergence, we still can
compare the performance of various algorithms. This has indeed formed a majority of
current research in algorithm development in the research community of optimization
and machine intelligence \cite{Enge,Talbi,Yang2}.

In combinatorial optimization, many important developments exist on complexity analysis, run time
and convergence analysis \cite{Reben,Neumann}. For continuous optimization, no-free-lunch-theorems do not hold \cite{Auger,Auger2}.
As a relatively young field, many open problems still remain in the field of randomized search heuristics \cite{Auger}.
In practice, most assume that metaheuristic algorithms tend to be
less complex for implementation, and in many cases, problem sizes are not directly linked with
the algorithm complexity. However, metaheuristics can often solve very tough NP-hard optimization,
while our understanding of the efficiency and convergence of metaheuristics lacks far behind.

Apart from the complex interactions among multiple search agents (making the mathematical
analysis intractable), another important issue is the various randomization techniques used for modern
metaheuristics, from simple randomization such as uniform distribution
to random walks, and to more elaborate L\'evy flights \cite{Blum,Pav,Yang2}. There
is no unified approach to analyze these mathematically.
In this paper, we intend to review the convergence of two metaheuristic algorithms
including simulated annealing and PSO, followed by the new convergence analysis of the
firefly algorithm. Then, we try to formulate
a framework for algorithm analysis in terms of Markov chain Monte Carlo.
We also try to analyze the mathematical and statistical foundations for randomization
techniques from simple random walks to L\'evy flights.
Finally, we will discuss some of important open questions as
further research topics.

\section{Convergence Analysis of Metaheuristics}

The formulation and numerical studies of various metaheuristics have been the main
focus of most research studies. Many successful applications have demonstrated the
efficiency of metaheuristics in various context, either through comparison with
other algorithms and/or applications to well-known problems. In contrast,
the mathematical analysis lacks behind, and convergence analysis has been
carried out for only a minority few algorithms such as simulated annealing and
particle swarm optimization \cite{Clerc,Neumann}. The main approach is often for very simplified systems
using dynamical theory and other ad hoc approaches. Here in this section, we first review the
simulated annealing and its convergence, and we move onto the population-based
algorithms such as PSO. We then take the recently developed
firefly algorithm as a further example to carry out its convergence analysis.

\subsection{Simulated Annealing}

Simulated annealing (SA) is one of the widely used metaheuristics, and is also one of
the most studies in terms of convergence analysis \cite{Berts,Kirk}. The essence
of simulated annealing is a trajectory-based random walk of a single agent,
starting from an initial guess $\x_0$. The next move only depends on the current
state or location and the acceptance probability $p$. This is essentially
a Markov chain whose transition probability from the current state to the next state
is given by
\be p=\exp \Big[ - \frac{\Delta E}{k_B T} \Big], \ee
where $k_B$ is Boltzmann's constant, and $T$ is the temperature.
Here the energy change $\Delta E$ can be linked with the change of
objective values. A few studies on the convergence of simulated annealing
have paved the way for analysis for all simulated annealing-based algorithms
\cite{Berts,Gran,Stein}. Bertsimas and Tsitsiklis provided an excellent review of
the convergence of SA under various assumptions \cite{Berts,Gran}.  By using the assumptions
that SA forms an inhomogeneous Markov chain with finite states, they proved a probabilistic convergence function $P$, rather than almost sure convergence, that
\be \max P \big[ \x_i(t) \in S_* \Big| \x_0\big] \ge \frac{A}{t^{\alpha}}, \ee
where $S_*$ is the optimal set, and $A$ and $\alpha$ are positive constants \cite{Berts}. This is for the cooling schedule $T(t) =d/\ln (t)$, where $t$ is the iteration counter
or pseudo time. These studies largely used Markov chains as the main tool.
We will come back later to a more general framework of Markov chain Monte Carlo (MCMC)
in this paper \cite{Gamer,Ghate}.

\subsection{PSO and Convergence}

Particle swarm optimization (PSO) was developed by Kennedy and
Eberhart in 1995 \cite{Kennedy,Kennedy2}, based on the swarm behaviour such
as fish and bird schooling in nature. Since then, PSO has
generated much wider interests, and forms an exciting, ever-expanding
research subject, called swarm intelligence. PSO has been applied
to almost every area in optimization, computational intelligence,
and design/scheduling applications.

The movement of a swarming particle consists of two major components:
a stochastic component and a deterministic component.
Each  particle is attracted toward the position of the current global best
$\ff{g}^*$ and its own best location $\x_i^*$ in history,
while at the same time it has a tendency to move randomly.

Let $\x_i$ and $\ff{v}_i$ be the position vector and velocity for
particle $i$, respectively. The new velocity and location updating formulas are
determined by
\be \ff{v}_i^{t+1}= \ff{v}_i^t  + \a \ff{\e}_1
[\ff{g}^*-\x_i^t] + \b \ff{\e}_2 [\x_i^*-\x_i^t].
\label{pso-speed-100}
\ee
\be \x_i^{t+1}=\x_i^t + \ff{v}_i^{t+1}, \label{pso-speed-140} \ee
where $\ff{\e}_1$ and $\ff{\e}_2$ are two random vectors, and each
entry taking the values between 0 and 1. The parameters $\a$ and $\b$ are the learning parameters or
acceleration constants, which can typically be taken as, say, $\a\=\b \=2$.

There are at least two dozen PSO variants which extend the standard PSO
algorithm, and the most noticeable improvement
is probably to use inertia function $\theta
(t)$ so that $\ff{v}_i^t$ is replaced by $\theta(t) \ff{v}_i^t$
where $\theta \in [0,1]$ \cite{Chat}.
This is equivalent to introducing a virtual mass to stabilize the motion
of the particles, and thus the algorithm is expected to converge more quickly.

The first convergence analysis of PSO was carried out by Clerc and Kennedy
in 2002 \cite{Clerc} using the theory of dynamical systems.
Mathematically, if we ignore the random factors, we can view the system
formed by (\ref{pso-speed-100}) and (\ref{pso-speed-140}) as a dynamical system.
If we focus on a single particle $i$ and imagine that there is only one particle in this system,
then the global best $\ff{g}^*$ is the same as its current best $\x_i^*$. In this case, we have
\be \ff{v}_i^{t+1} = \ff{v}_i^t + \gamma (\ff{g}^*-\x_i^t), \quad \gamma=\a+\b, \ee
and
\be \x_i^{t+1} =\x_i^t + \ff{v}_i^{t+1}. \ee
Considering the 1D dynamical system for particle swarm optimization,
we can replace $\ff{g}^*$ by a parameter constant $p$ so
that we can see if or not the particle of interest will converge towards $p$.
By setting $u_t=p-x(t+1)$
and using the notations for dynamical systems, we have a simple dynamical system
\be v_{t+1}=v_t + \gamma u_t, \quad u_{t+1} =-v_t + (1-\gamma) u_t, \ee
or
\be Y_{t+1} = A Y_t, \quad A=\mat{ 1 & & \gamma \\ -1 & & 1-\gamma}, \quad Y_t=\mat{v_t \\ u_t}. \ee
The general solution of this dynamical system can be written as $Y_t=Y_0 \exp[A t]$.
The system behaviour can be characterized by the eigenvalues $\lam$ of $A$, and we have
$\lam_{1,2}=1-\gamma/2 \pm \sqrt{\gamma^2 - 4 \gamma}/2$.
It can be seen clearly that $\gamma=4$ leads to a  bifurcation.
Following a straightforward analysis of this dynamical system,
we can have three cases. For $0 < \gamma <4$, cyclic and/or quasi-cyclic
trajectories exist. In this case, when randomness is gradually reduced,
some convergence can be observed.
For $\gamma>4$, non-cyclic behaviour can be expected and
the distance from $Y_t$ to the center $(0,0)$ is monotonically increasing with $t$.
In a special case $\gamma=4$, some convergence behaviour can be observed.
For detailed analysis, please refer to Clerc and Kennedy \cite{Clerc}.
Since $p$ is linked with the global best, as the iterations continue,
it can be expected that all particles will aggregate towards the
the global best.

\subsection{Firefly algorithm, Convergence and Chaos}

Firefly Algorithm (FA) was developed by Yang \cite{Yang,Yang3},
which was based on the flashing patterns and behaviour
of fireflies. In essence, each firefly will be attracted to brighter ones,
while at the same time, it explores and searches for prey randomly.
In addition, the brightness of a firefly is determined by the landscape of the
objective function.

The movement of a firefly $i$ is attracted to another more attractive (brighter)
firefly $j$ is determined by
\be \x_i^{t+1} =\x_i^t + \b_0 e^{-\gamma r^2_{ij}} (\x_j^t-\x_i^t) + \a \; \ff{\e}_i^t,
\label{FA-equ-50}  \ee
where the second term is due to the attraction. The third term
is randomization with $\a$ being the randomization parameter, and
$\ff{\e}_i^t$ is a vector of random numbers drawn from a Gaussian distribution
or other distributions.  Obviously, for a given firefly, there are often many more
attractive fireflies, then we can either go through all of them via a loop or
use the most attractive one. For multiple modal problems, using a loop while moving toward
each brighter one is usually more effective, though this will lead to a slight increase
of algorithm complexity.

Here is $\b_0 \in [0,1]$ is the attractiveness at $r=0$,
and $r_{ij}=||\x_i-\x_j||_2$ is the $\ell_2$-norm or Cartesian distance. For other problems such as scheduling,
any measure that can effectively characterize the quantities of interest in the optimization
problem can be used as the `distance' $r$.

For most implementations, we can take $\b_0=1$, $\a=O(1)$ and $\gamma=O(1)$.
It is worth pointing out that (\ref{FA-equ-50}) is essentially a random walk biased
towards the brighter fireflies. If $\b_0=0$, it becomes a simple random walk.
Furthermore, the randomization term can easily be
extended to other distributions such as L\'evy flights \cite{Gutow,Pav}.

\begin{figure}
\centerline{\includegraphics[height=2in,width=2.25in]{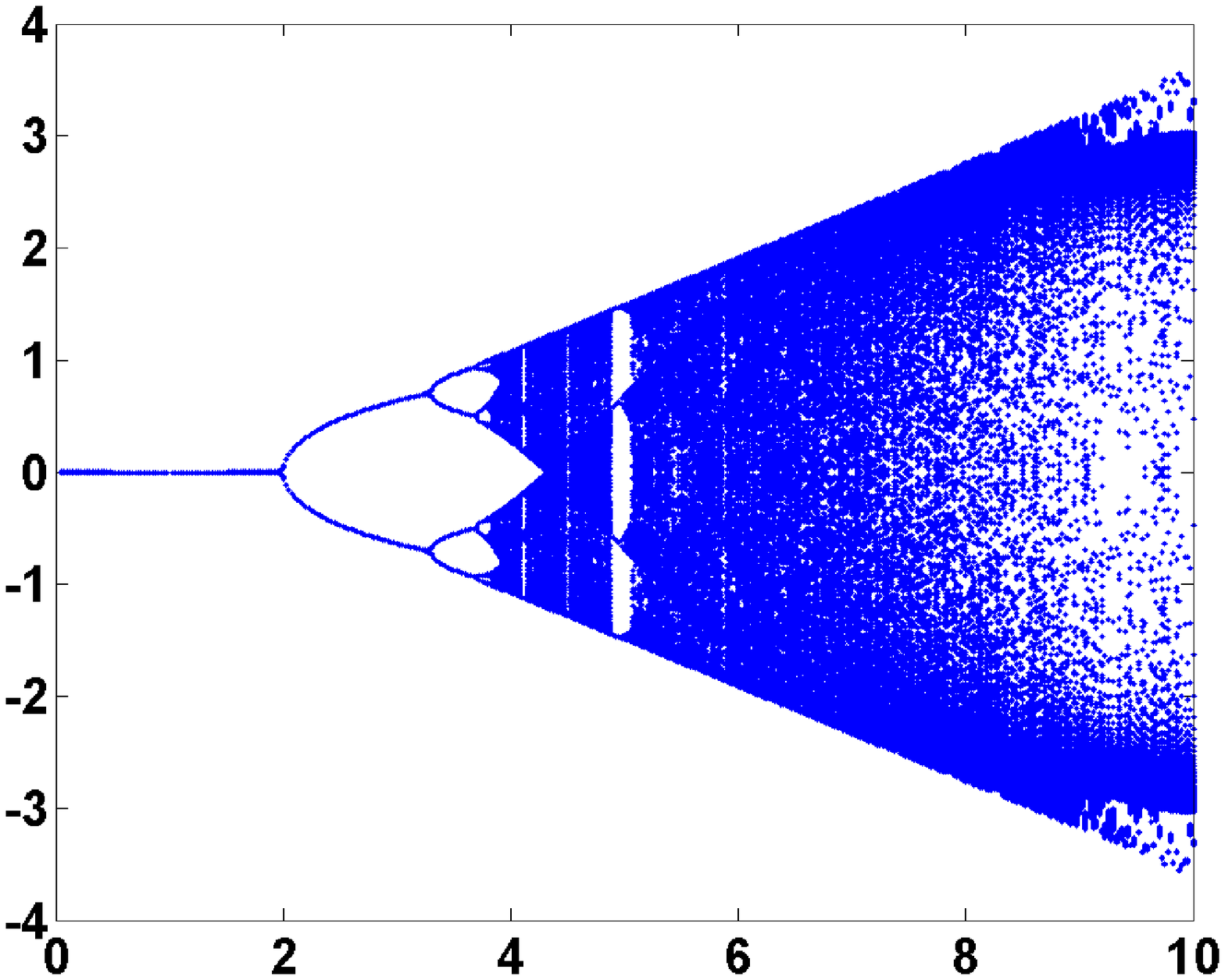}
\includegraphics[height=2in,width=2.5in]{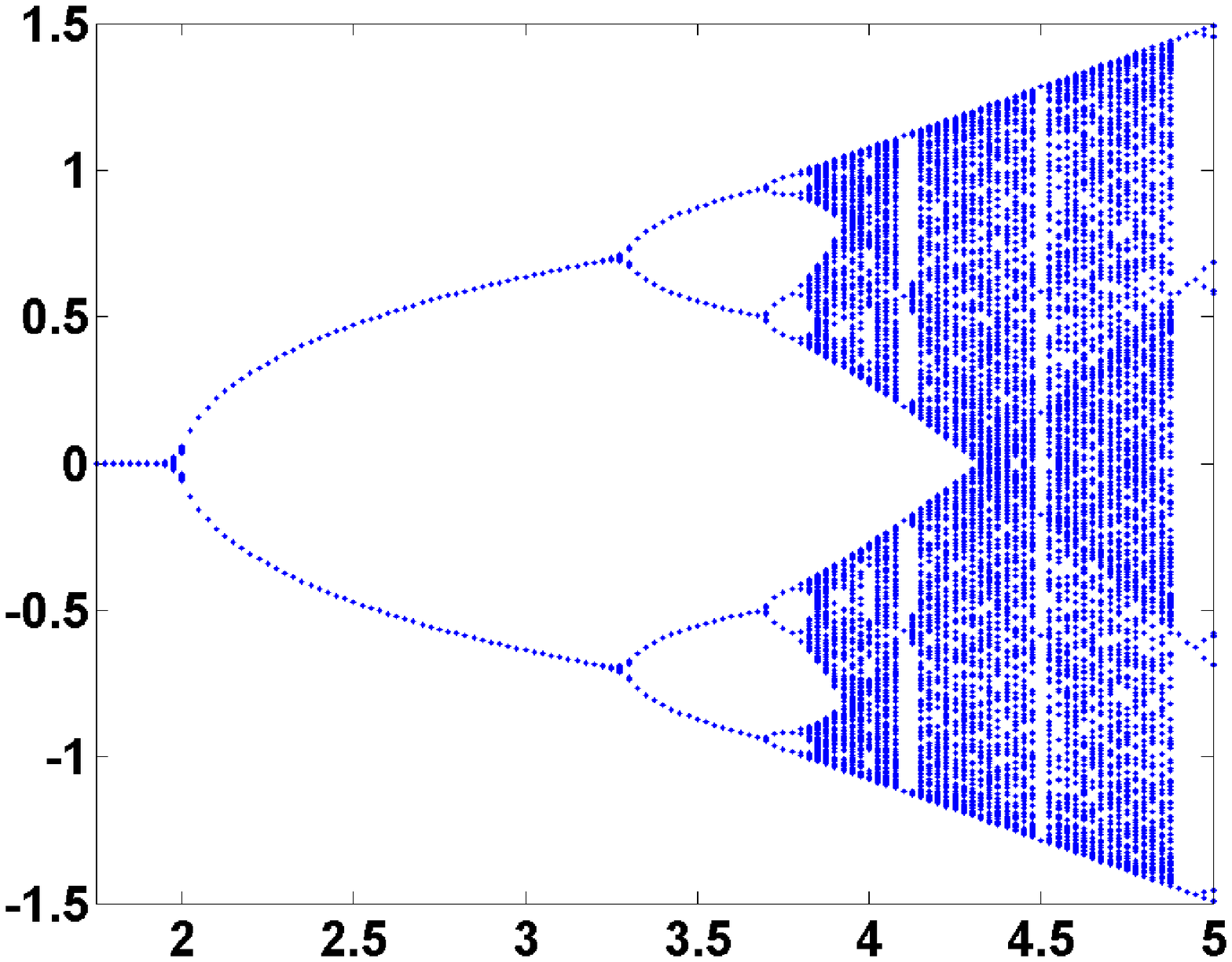}}
\caption{The chaotic map of the iteration formula (\ref{fa-dynamics-150}) in the firefly algorithm
and the transition between from periodic/multiple  states to chaos. \label{fap-fig-100} }
\end{figure}

We now can carry out the convergence analysis for the firefly algorithm in a
framework similar to Clerc and Kennedy's dynamical analysis. For simplicity,
we start from the equation for firefly motion without the randomness term
\be \x_i^{t+1}=\x_i^t + \b_0 e^{-\gamma r_{ij}^2} (\x_j^t - \x_i^t). \ee
If we focus on a single agent, we can replace $\x_j^t$ by the global best
$g$ found so far, and we have
\be \x_i^{t+1} =\x_i^t + \b_0 e^{-\gamma r_{i}^2 } (g-\x_i^t), \ee
where the distance $r_i$ can be given by the $\ell_2$-norm $r_i^2=||g-\x_i^t||_2^2$.
In an even simpler 1-D case, we can set $y_{t}=g-\x_i^t$, and we have
\be y_{t+1} =y_t - \b_0 e^{-\gamma y_t^2} y_t. \label{fa-dynamics-25} \ee
We can see that $\gamma$ is a scaling parameter which only affects the scales/size
of the firefly movement. In fact, we can let
$u_t=\sqrt{\gamma} y_t$ and we have
\be u_{t+1}=u_t [1-\b_0 e^{-u_t^2}]. \label{fa-dynamics-150} \ee
These equations can be analyzed easily using the same methodology for studying
the well-known logistic map \be u_{t+1}=\lam u_t (1-u_t). \label{fa-chaos-55} \ee
The chaotic map is shown
in Fig. \ref{fap-fig-100}, and the focus on the transition from periodic multiple
states to chaotic behaviour is shown in the same figure.

As we can see from Fig. \ref{fap-fig-100} that convergence can be achieved
for $\b_0 <2$. There is a transition from periodic to chaos at $\b_0 \approx 4$.
This may be surprising, as the aim of designing a metaheuristic algorithm is
to try to find the optimal solution efficiently and accurately. However,
chaotic behaviour is not necessarily a nuisance; in fact, we can use it to
the advantage of the firefly algorithm. Simple chaotic characteristics from  (\ref{fa-chaos-55})
can often be used as an efficient mixing technique for generating diverse solutions.
Statistically, the logistic mapping (\ref{fa-chaos-55}) with $\lam=4$ for the initial
states in (0,1) corresponds a beta distribution
\be B(u,p,q)=\frac{\Gamma(p+q)}{\Gamma(p) \Gamma(q)} u^{p-1} (1-u)^{q-1}, \ee
when $p=q=1/2$. Here $\Gamma(z)$ is the Gamma function
\be \Gamma(z) =\int_0^{\infty} t^{z-1} e^{-t} dt.  \ee
In the case when  $z=n$ is an integer, we have $\Gamma(n)=(n-1)!$. In addition,
$\Gamma(1/2)=\sqrt{\pi}$. From the algorithm implementation point of view,
we can use higher attractiveness $\b_0$ during the early stage
of iterations so that the fireflies can explore, even chaotically,
the search space more effectively. As the search continues and convergence
approaches, we can reduce the attractiveness $\b_0$ gradually, which may
increase the overall efficiency of the algorithm. Obviously, more
studies are highly needed to confirm this.

\subsection{Markov Chain Monte Carlo}

From the above convergence analysis, we know that
there is no mathematical framework in general to provide insights into the
working mechanisms, the stability and convergence of a give algorithm.
Despite the increasing popularity of metaheuristics, mathematical analysis
remains fragmental, and many open problems need urgent attention.

Monte Carlo methods have been applied in many applications \cite{Sobol},
including almost all areas of sciences and engineering. For example, Monte
Carlo methods are widely used in uncertainty and sensitivity analysis \cite{Matthews}.
From the statistical point of view, most metaheuristic
algorithms can be viewed in the framework of Markov chains \cite{Ghate,Sobol}. For example,
simulated annealing \cite{Kirk} is a Markov
chain, as the next state or new solution
in SA only depends on the current state/solution and the transition probability.
For a given Markov chain with certain ergodicity, a stability probability distribution
and convergence can be achieved.

Now if look at the PSO closely using the framework of Markov chain
Monte Carlo \cite{Gamer,Geyer,Ghate}, each particle in PSO essentially
forms a Markov chain, though this Markov chain is biased towards to the
current best, as the transition probability often leads to the
acceptance of the move towards the current global best.
Other population-based algorithms can also be viewed in this framework.
In essence, all metaheuristic algorithms with piecewise, interacting paths
can be analyzed in the general framework of Markov chain Monte Carlo.
The main challenge is to realize this and to use the appropriate Markov chain
theory to study metaheuristic algorithms. More fruitful studies will surely emerge in the future.

\section{Search Efficiency and Randomization}

Metaheuristics can be considered as an efficient way to produce acceptable
solutions by trial and error to a complex problem in
a reasonably practical time. The complexity of the problem of interest makes it impossible
to search every possible solution or combination, the aim is
to find good feasible solutions in an acceptable timescale. There is
no guarantee that the best solutions can be found, and we even do not
know whether an algorithm will work and why if it does work. The idea is
to have an efficient but practical algorithm that will work most
the time and is able to produce good quality solutions. Among the
found quality solutions, it is expected some of them are nearly optimal,
though there is no guarantee for such optimality.

The main components of any metaheuristic algorithms are: intensification
and diversification, or exploitation and exploration \cite{Blum,Yang2}.
Diversification means to generate diverse solutions so as to explore the search space
on the global scale, while intensification means to focus on the
search in a local region by exploiting the information that a current
good solution is found in this region. This is in combination with
with the selection of the best solutions.

As discussed earlier, an important component in swarm intelligence
and modern metaheuristics is randomization, which enables an algorithm
to have the ability to jump out of any local optimum so as to search globally.
Randomization can also be used for local search around the current best
if steps are limited to a local region. Fine-tuning the randomness
and balance of local search and global search is crucially important
in controlling the performance of any metaheuristic algorithm.

Randomization techniques can be a very simple method using uniform distributions,
or more complex methods as those used in Monte Carlo simulations \cite{Sobol}. They can also
be more elaborate, from Brownian random walks to L\'evy flights.

\subsection{Gaussian Random Walks}

A random walk is a random process which consists of taking a series of consecutive
random steps. Mathematically speaking, let $u_N$ denotes the sum of each
consecutive random step $s_i$, then $u_N$ forms a random walk
\be u_N=\sum_{i=1}^N s_i = s_1 + ... + s_N = u_{N-1} + s_N, \label{Walk-markov-50}  \ee
where $s_i$ is a random step drawn from a random distribution. This suggests that
the next state $u_N$ will only depend the current existing state $u_{N-1}$
and the motion or transition $u_N$ from the existing state to the next state.
In theory, as the number of steps $N$ increases, the central limit theorem implies that
the random walk (\ref{Walk-markov-50}) should approaches a Gaussian distribution.
In addition, there is no reason why each step length
should be fixed. In fact, the step size can also vary according to a known distribution.
If the step length obeys the Gaussian distribution, the random walk becomes the standard
Brownian motion \cite{Gutow,Yang2}.

From metaheuristic point of view, all paths of search agents form a random walk,
including a particle's trajectory in simulated annealing,
a zig-zag path of a particle in PSO, or the piecewise path of a firefly in FA.
The only difference is that transition probabilities are different, and change
with time and locations.

Under simplest assumptions, we know that a Gaussian distribution is stable.
For a particle starts with an initial location $\x_0$, its final location $\x_N$
after $N$ time steps is
\be \x_N =\x_0 +  \sum_{i=1}^N \alpha_i s_i, \ee
where $\alpha_i >0$ is a parameters controlling the step sizes or scalings.
If $s_i$ is drawn from a normal distribution $N(\mu_i, \sigma_i^2)$,
then the conditions of stable distributions lead to a combined Gaussian distribution
\be \x_N \sim N(\mu_*, \sigma_*^2), \quad
\mu_* = \sum_{i=1}^N \alpha_i \mu_i, \quad \sigma_*^2=\sum_{i=1}^N \alpha_i [\sigma_i^2
+(\mu_*-\mu_i)^2]. \ee We can see that that the mean location changes with $N$ and
the variances increases as $N$ increases, this makes it possible to reach any areas in
the search space if $N$ is large enough.

A diffusion process can be viewed as a series of Brownian motion, and the motion
obeys the Gaussian distribution. For this reason, standard diffusion is often referred to as the Gaussian diffusion. As the mean of particle locations
is obviously zero if $\mu_i=0$, their variance will increase linearly with $t$.
In general, in the $d$-dimensional space, the variance of Brownian random walks can be written as
\be \sigma^2(t) = |v_0|^2 t^2 + (2 d D) t, \ee
where $v_0$ is the drift velocity of the system. Here $D=s^2/(2 \tau)$ is the effective diffusion coefficient
which is related to the step length $s$ over a short time interval $\tau$ during each jump.
If the motion at each step is not Gaussian, then the diffusion is called non-Gaussian diffusion.
If the step length obeys other distribution, we have to deal with more generalized random walks. A very special case
is when the step length obeys the L\'evy distribution, such a random walk is called  L\'evy flight or L\'evy walk.

\subsection{Randomization via L\'evy Flights}

In nature, animals search for food in a random or quasi-random manner.
In general, the foraging path of an animal is effectively a random walk because the next move is
based on the current location/state and the transition probability to
the next location. Which direction it chooses depends implicitly on a
probability which can be modelled mathematically \cite{Bell,Pav}.
For example, various studies
have shown that the flight behaviour of many animals and insects has demonstrated
the typical characteristics of L\'evy flights \cite{Pav,Reynolds}.
Subsequently, such behaviour has been applied to
optimization and optimal search, and preliminary results show its
promising capability \cite{Viswan,Pav}.

In general, L\'evy distribution is stable, and can be defined in terms of a characteristic function or
the following Fourier transform
\be F(k)=\exp[-\a |k|^{\b}], \quad 0 < \b \le 2, \ee
where $\a$ is a scale parameter. The inverse of this integral is not easy, as it does not
have nay analytical form, except for a few special cases \cite{Gutow,Nol}.
For the case of $\b=2$, we have $F(k)=\exp[-\a k^2]$,
whose inverse Fourier transform corresponds to a Gaussian distribution.
Another special case is $\b=1$, which corresponds to a Cauchy distribution

For the general case, the inverse integral
\be L(s) =\frac{1}{\pi} \int_0^{\infty} \cos (k s) \exp[-\a |k|^{\b}] dk, \ee
can be estimated only when $s$ is large. We have
\be L(s) \rightarrow \frac{\a \; \b \; \Gamma(\b) \sin (\pi \b/2)}{ \pi |s|^{1+\b}}, \quad s \rightarrow \infty. \ee

L\'evy flights are more efficient than Brownian random walks
in exploring unknown, large-scale search space.
There are many reasons to explain this efficiency, and one of them
is due to the fact that the variance of L\'evy flights takes the following form
\be \sigma^2(t) \sim t^{3-\b},  \quad 1 \le \b \le 2, \ee
which increases much faster than the linear relationship (i.e., $\sigma^2(t) \sim t$) of Brownian random walks.

Studies show that L\'evy flights can maximize the efficiency of resource searches in uncertain environments.
In fact,  L\'evy flights have been observed among
foraging patterns of albatrosses and fruit flies \cite{Pav,Reynolds,Viswan}.  In addition,
L\'evy flights have many applications. Many physical phenomena such as
the diffusion of fluorenscent molecules, cooling behavior and noise could show L\'evy-flight
characteristics under the right conditions \cite{Reynolds}.

\section{Open Problems}

It is no exaggeration to say that metahueristic algorithms have been a great success
in solving various tough optimization problems. Despite this huge success,
there are many important questions which remain unanswered.
We know how these heuristic algorithms work, and we also partly
understand why these algorithms work. However,  it is difficult
to analyze mathematically why these algorithms are so successful,
though significant progress has been made in the last few years \cite{Auger,Neumann}.
However, many open problems still remain.

For all population-based metaheuristics, multiple search agents form
multiple interacting Markov chains. At the moment, theoretical development in these areas
are still at early stage. Therefore, the mathematical analysis concerning of the
rate of convergence is very difficult, if not impossible.
Apart from the mathematical analysis on a limited few metaheuristics,
convergence of all other algorithms has not been proved mathematically, at least
up to now. Any mathematical analysis will thus provide important insight into these algorithms. It will also be
valuable for providing new directions for further important modifications on these algorithms
or even pointing out innovative ways of developing new algorithms.

For almost all metaheuristics including future new algorithms, an important issue
to be addresses is to provide a balanced trade-off between local intensification
and global diversification \cite{Blum}. At present, different algorithm uses different
techniques and mechanism with various parameters to control this, they are far from
optimal.  Important questions are: Is there any optimal way to achieve this balance?
If yes, how? If not, what is the best we can achieve?

Furthermore, it is still only partly understood why different components of heuristics and metaheuristics
interact in a coherent and balanced way so that they produce efficient algorithms which converge under the given
conditions. For example, why does a balanced combination of randomization and a deterministic component
lead to a much more efficient algorithm (than a purely deterministic and/or a purely random algorithm)?
How to measure or test if a balance is reached? How to prove that the use of memory can significantly increase the
search efficiency of an algorithm? Under what conditions?

In addition, from the well-known No-Free-Lunch theorems \cite{Wolp},
we know that they have been proved for single objective optimization for finite search domains,
but they do not hold for continuous infinite domains \cite{Auger,Auger2}.
In addition, they remain unproved for multiobjective optimization. If they are proved to be true (or not)
for multiobjective optimization, what are the implications for algorithm development?
Another important question is about the performance comparison. At the moment,
there is no agreed measure for comparing performance of different algorithms, though
the absolute objective value and the number of function evaluations are two widely used
measures.  However, a formal theoretical analysis is yet to be developed.

Nature provides almost unlimited ways for problem-solving. If we can observe carefully,
we are surely inspired to develop more powerful and efficient new generation algorithms.
Intelligence is a product of biological evolution in nature.
Ultimately some intelligent algorithms (or systems) may appear in the future,
so that they can evolve and optimally adapt to solve
NP-hard optimization problems efficiently and intelligently.

Finally, a current trend is to use simplified metaheuristic algorithms to
deal with complex optimization problems. Possibly, there is a need to
develop more complex metaheuristic algorithms which can truly mimic
the exact working mechanism of some natural or biological systems,
leading to more powerful next-generation, self-regulating, self-evolving,
and truly intelligent metaheuristics.

\end{document}